\documentclass[11pt]{amsart}
\RequirePackage{url}
\usepackage{verbatim}
\usepackage {   amsfonts    }
\usepackage {   amsmath     }
\usepackage {   amsthm      }
\usepackage {   latexsym    }
\usepackage {   graphics    }
\usepackage {   color       }
\usepackage {   multicol    }
\usepackage {   amssymb     }
\usepackage {   enumerate   }

\makeatletter

\@namedef{subjclassname@2010}{%

  \textup{2010} Mathematics Subject Classification}

\makeatother

\newtheorem{thm}{Theorem}[section]

\theoremstyle{definition}

\numberwithin{equation}{section}

\newcommand{\prf}[1]{\noindent \emph{Proof.} #1 \qed\bigskip}

\begin{document}

\baselineskip=17pt

\title[4-Colored F-Partitions Mod 5]{An Unexpected Congruence Modulo 5 for 4--Colored Generalized Frobenius Partitions}

\author[J. A. Sellers]{James A. Sellers}
\address{Department of Mathematics, Penn State University, University Park, PA  16802, USA, sellersj@psu.edu}
\thanks{J. A. Sellers gratefully acknowledges the support of the Austrian American Educational Commission which supported him during the Summer Semester 2012 as a Fulbright Fellow at the Johannes Kepler University, Linz, Austria.}

\dedicatory{On the occasion of the $125^{th}$ anniversary of the birth of Srinivasa Ramanujan}

\date{\today}

\begin{abstract}
In his 1984 AMS Memoir, George Andrews defined the family of $k$--colored generalized Frobenius partition functions.  These are denoted by  $c\phi_k(n)$ where $k\geq 1$ is the number of colors in question.  In that Memoir, Andrews proved (among many other things) that, for all $n\geq 0,$ $c\phi_2(5n+3) \equiv 0\pmod{5}.$  Soon after, many authors proved congruence properties for various $k$--colored generalized Frobenius partition functions, typically with a small number of colors.  

In 2011, Baruah and Sarmah proved a number of congruence properties for $c\phi_4$, all with moduli which are powers of 4.  In this brief note, we add to the collection of congruences for $c\phi_4$ by proving this function satisfies an unexpected result modulo 5.  The proof is elementary, relying on Baruah and Sarmah's results as well as work of Srinivasa Ramanujan.  
\end{abstract}

\maketitle
\section{Introduction}
In his 1984 AMS Memoir, George Andrews \cite{AndMem} defined the family of $k$--colored generalized Frobenius partition functions which are denoted by  $c\phi_k(n)$ where $k\geq 1$ is the number of colors in question.  Among many things, 
Andrews \cite[Corollary 10.1]{AndMem} proved that, for all $n\geq 0,$ $c\phi_2(5n+3) \equiv 0\pmod{5}.$  Soon after, many authors proved similar congruence properties for various $k$--colored generalized Frobenius partition functions, typically for a small number of colors $k.$  See, for example, \cite{ES, GarThesis, Kol1, KolPow3, Lovejoy, Ono, PauRad, Sel1, Xiong}.

In 2011, Baruah and Sarmah \cite{BarSar} proved a number of congruence properties for $c\phi_4$, all with moduli which are powers of 4.  In this brief note, we add to the collection of congruences for $c\phi_4$ by proving this function satisfies an unexpected result modulo 5.  

\begin{thm}
\label{congmod5}
For all $n\geq 0,$ $c\phi_4(10n+6) \equiv 0 \pmod{5}.$
\end{thm}
Our proof is elementary, relying on Baruah and Sarmah's results as well as work of Srinivasa Ramanujan. 

\section{An Elementary Proof of Theorem \ref{congmod5}}
Recall Ramanujan's functions 
$$
\phi(q) := \sum_{n=-\infty}^\infty q^{n^2}  \text{\ \ and \ \ } 
\psi(q) := \sum_{n=0}^\infty q^{n(n+1)/2}.  
$$
Using Jacobi's Triple Product Identity \cite[Entry 19]{Berndt}, we have the following well--known product representations for $\phi(q)$ and $\psi(q):$  
\begin{equation}
\label{phi_prod}
\phi(q) = \frac{(q^2;q^2)_\infty^5}{(q;q)_\infty^2(q^4;q^4)_\infty^2}
\end{equation}
and
\begin{equation}
\label{psi_prod}
\psi(q) = \frac{(q^2;q^2)_\infty^2}{(q;q)_\infty}
\end{equation}

Baruah and Sarmah \cite[Theorem 2.1]{BarSar} proved the following valuable representation of the generating function for $c\phi_4.$  
\begin{thm}
\label{genfncphi4}
$$
\sum_{n=0}^\infty c\phi_4(n)q^n = \frac{\phi^3(q^2) + 12q\phi(q^2)\psi^2(q^4)}{(q;q)_\infty^4}
$$
where $(a;b)_\infty := (1-a)(1-ab)(1-ab^2)(1-ab^3)\dots$  
\end{thm}
From here, we wish to 2--dissect the generating function in Theorem \ref{genfncphi4} (because we want to study the coefficients of $q^{10n+6}$ in the power series representation of the generating function for $c\phi_4(n)$).  To complete this task, we follow the path laid out by Baruah and Sarmah \cite{BarSar}.  We begin by rewriting the generating function in Theorem \ref{genfncphi4} as
$$
\sum_{n=0}^\infty c\phi_4(n)q^n = \frac{\phi^3(q^2) + 12q\phi(q^2)\psi^2(q^4)}{(q;q^2)_\infty^4(q^2;q^2)_\infty^4}.  
$$
Then we see that 
\begin{eqnarray*}
&& 
\sum_{n=0}^\infty c\phi_4(n)q^n + \sum_{n=0}^\infty c\phi_4(n)(-q)^n \nonumber \\
&=& 
\frac{\phi^3(q^2)}{(q^2;q^2)_\infty^4}\left\{ \frac{1}{(q;q^2)_\infty^4} + \frac{1}{(-q;q^2)_\infty^4} \right\} 
+
12q\frac{\phi(q^2)\psi^2(q^4)}{(q^2;q^2)_\infty^4}\left\{ \frac{1}{(q;q^2)_\infty^4} - \frac{1}{(-q;q^2)_\infty^4} \right\} \nonumber \\
&=& 
\frac{\phi^3(q^2)}{(q^2;q^2)_\infty^4(q^2;q^4)_\infty^4}\left\{ (-q;q^2)_\infty^4 +(q;q^2)_\infty^4 \right\} 
+
12q\frac{\phi(q^2)\psi^2(q^4)}{(q^2;q^2)_\infty^4(q^2;q^4)_\infty^4}\left\{ (-q;q^2)_\infty^4 - (q;q^2)_\infty^4 \right\}. 
\end{eqnarray*}
As noted by Baruah and Sarmah \cite[(3.11) and (3.12)]{BarSar}, we can employ work of Ramanujan \cite[Entry 25]{Berndt} to obtain
\begin{equation*}
(-q;q^2)_\infty^4 +(q;q^2)_\infty^4 = 2\frac{\phi^2(q^2)}{(q^2;q^2)_\infty^2}
\end{equation*}
and
\begin{equation*}
(-q;q^2)_\infty^4 -(q;q^2)_\infty^4 = 8q\frac{\psi^2(q^4)}{(q^2;q^2)_\infty^2}.
\end{equation*}
These can be used in the above to obtain, after simplification, 
$$
\sum_{n=0}^\infty c\phi_4(n)q^n + \sum_{n=0}^\infty c\phi_4(n)(-q)^n 
=
2\left\{ \frac{\phi^5(q^2)}{(q^2;q^2)_\infty^6 (q^2;q^4)_\infty^4}+48q^2\frac{\phi(q^2)\psi^4(q^4)}{(q^2;q^2)_\infty^6 (q^2;q^4)_\infty^4} \right\}.
$$
Therefore, 
\begin{equation}
\label{2-dissection1}
\sum_{n\geq 0} c\phi_4(2n)q^n = \frac{\phi^5(q)}{(q;q)_\infty^6 (q;q^2)_\infty^4}+48q\frac{\phi(q)\psi^4(q^2)}{(q;q)_\infty^6 (q;q^2)_\infty^4}.
\end{equation}
We now utilize (\ref{phi_prod}) and (\ref{psi_prod}) in (\ref{2-dissection1}) to obtain 
\begin{equation}
\label{2-dissection2}
\sum_{n\geq 0} c\phi_4(2n)q^n = \frac{(q^2;q^2)^{29}}{(q;q)_\infty^{20} (q^4;q^4)_\infty^{10}}+48q\frac{(q^2;q^2)_\infty^5(q^4;q^4)_\infty^6}{(q;q)_\infty^{12}}
\end{equation}
after elementary simplifications. 

Our goal now is to prove that, when written as power series in $q,$ each of the two functions on the right--hand side of (\ref{2-dissection2}) has the property that every coefficient of a term of the form $q^{5n+3}$ is a multiple of 5.  Then Theorem \ref{congmod5} follows.  

With this in mind, we first note that  
\begin{equation}
\label{firstterm_mod5}
\frac{(q^2;q^2)^{29}}{(q;q)_\infty^{20} (q^4;q^4)_\infty^{10}}
\equiv 
\frac{(q^{10};q^{10})^{5}(q^2;q^2)_\infty^4}{(q^5;q^5)_\infty^{4} (q^{20};q^{20})_\infty^{2}}\pmod{5}.
\end{equation}
Now all of the functions of $q^5$ in the right--hand side of (\ref{firstterm_mod5}) can be ignored; in other words, if we can show that every coefficient of a term of the form $q^{5n+3}$ is a multiple of 5 in the power series representation of the function $(q^2;q^2)_\infty^4,$ then every coefficient of a term of the form $q^{5n+3}$ is a multiple of 5 in the power series representation of the function 
$$
\frac{(q^{10};q^{10})^{5}(q^2;q^2)_\infty^4}{(q^5;q^5)_\infty^{4} (q^{20};q^{20})_\infty^{2}}
$$
as well.  So we focus our attention on $(q^2;q^2)_\infty^4.$  We know from Euler's Pentagonal Number Theorem and a well--known result of Jacobi \cite[page 11 and page 176]{Andrews} that 
\begin{equation}
\label{PNTandGauss1}
(q^2;q^2)_\infty^4 = \sum_{m=-\infty}^\infty \sum_{k=0}^\infty (-1)^{k+m}(2k+1)q^{k(k+1)+m(3m-1)}.
\end{equation}
We now consider exponents of $q$ and want to know when 
\begin{equation}
\label{PNTandGauss2}
5n+3 = k(k+1)+m(3m-1)
\end{equation}
has a solution for integers $k, m,$ and $n.$  When we consider (\ref{PNTandGauss2}) mod 5, we have 
\begin{equation*}
3 \equiv k(k+1)+m(3m-1)\pmod{5}.  
\end{equation*}
Completing the square gives 
\begin{equation}
\label{PNTandGauss4}
0 \equiv 3(2k+1)^2+(6m-1)^2\pmod{5}.  
\end{equation}
We know that all squares are congruent to 0, 1, or 4 modulo 5.  So we consider the nine possible cases in (\ref{PNTandGauss4}) and note rather quickly that $(2k+1)^2$ must be congruent to 0 modulo 5 (the other possibilities lead to no solutions).  Thus, $2k+1$ is divisible by 5 whenever we consider the coefficients of the terms of the form $q^{5n+3}.$  Therefore, by (\ref{PNTandGauss1}), we know that the coefficients of the terms of the form $q^{5n+3}$ in (\ref{firstterm_mod5}) must be divisible by 5.  

We turn our attention to the second term on the right--hand side of (\ref{2-dissection2}) and note that 
\begin{eqnarray}
\label{secondterm_mod5}
48q\frac{(q^2;q^2)_\infty^5(q^4;q^4)_\infty^6}{(q;q)_\infty^{12}}
&=&
48q\frac{(q^2;q^2)_\infty^5(q;q)_\infty^3(q^4;q^4)_\infty^6}{(q;q)_\infty^{15}}\nonumber \\
&\equiv &
48q\frac{(q^{10};q^{10})(q^{20};q^{20})_\infty(q;q)_\infty^3(q^4;q^4)_\infty}{(q^{5};q^{5})_\infty^{3}}\pmod{5}.
\end{eqnarray}
As before, we can ignore those functions which are functions of $q^5$, so our goal is to focus on the coefficients of the terms of the form $q^{5n+3}$ in 
$$
q(q;q)_\infty^3(q^4;q^4)_\infty
$$
which is equivalent to considering the coefficients of the terms of the form $q^{5n+2}$ in 
$$
(q;q)_\infty^3(q^4;q^4)_\infty.
$$
As above, we know that 
\begin{equation}
\label{PNTandGaussAgain}
(q;q)_\infty^3(q^4;q^4)_\infty = \sum_{m=-\infty}^\infty \sum_{k=0}^\infty (-1)^{k+m}(2k+1)q^{k(k+1)/2+2m(3m-1)}.
\end{equation}
Thus, we want to know when 
\begin{equation}
\label{PNTandGauss6}
5n+2 = \frac{k(k+1)}{2}+2m(3m-1)
\end{equation}
has a solution for integers $k, m,$ and $n.$ 
When we consider (\ref{PNTandGauss6}) modulo 5, we have 
\begin{equation}
\label{PNTandGauss7}
2 = \frac{k(k+1)}{2}+2m(3m-1) \pmod{5}
\end{equation}
and completing the square in (\ref{PNTandGauss7}) gives 
\begin{equation}
\label{PNTandGauss8}
0 = (2k+1)^2+8(m-1)^2 \pmod{5}.
\end{equation}
As above, we find that the only way that (\ref{PNTandGauss8}) can have solutions is if $2k+1\equiv 0 \pmod{5}.$  Thus, thanks to (\ref{PNTandGaussAgain}), we can conclude that all the coefficients of the terms of the form $q^{5n+3}$ in 
$$
48q\frac{(q^2;q^2)_\infty^5(q^4;q^4)_\infty^6}{(q;q)_\infty^{12}}
$$
are also divisible by 5.  This completes the proof of Theorem \ref{congmod5}.  
\hfill\qed

\section{Concluding Thoughts}
We close with a number of thoughts.  First, it is worth noting that Theorem \ref{congmod5} implies another congruence property modulo 5, this time for the function $\phi_4(n)$ (which is the number of generalized Frobenius partitions of $n$ which allow up to 4 repetitions of an integer in either row).  See Andrews \cite{AndMem} for more details.  We can easily prove the following result: 

\begin{thm}
\label{phi4mod5}
For all $n\geq 0,$ $\phi_4(10n+6) \equiv 0 \pmod{5}.$
\end{thm}

\prf{
Garvan \cite{GarThesis} proved that, for any prime $p,$ 
$$
\phi_{p-1}(n) \equiv c\phi_{p-1}(n) \pmod{p}
$$
for any integer $n\geq 0.$  Thus, Theorem \ref{phi4mod5} immediately follows from Theorem \ref{congmod5}. 
}

Secondly, we note that exactly the same kind of proof as that given for Theorem \ref{congmod5} can be employed to prove the following theorem:
\begin{thm}
\label{cphibar4mod5}
For all $n\geq 0,$ $\overline{c\phi}_4(10n+6) \equiv 0 \pmod{5}.$
\end{thm}
Note that the functions $\overline{c\phi}_k(n)$ were defined by Kolitsch \cite{Kol2} for any integer $k\geq 1.$  Indeed, Baruah and Sarmah \cite[(5.8)]{BarSar} prove that 
\begin{equation*}
\sum_{n=0}^\infty \overline{c\phi}_4(2n)q^n = 64q\frac{(q^4;q^4)_\infty^6}{(q^2;q^2)_\infty^7(q;q^2)_\infty^{12}}.  
\end{equation*}
As in the above work, we then see that 
\begin{eqnarray}
\label{cphibar2}
\sum_{n=0}^\infty \overline{c\phi}_4(2n)q^n 
&=& 
64q\frac{(q^4;q^4)_\infty^6(q^2;q^2)_\infty^5}{(q;q)_\infty^{12} } \nonumber \\
&\equiv &
64q\frac{(q^{20};q^{20})_\infty(q^4;q^4)_\infty(q^{10};q^{10})_\infty(q;q)_\infty^3}{(q^5;q^5)_\infty^{3}} \pmod{5}.
\end{eqnarray}
Thus, thanks to (\ref{cphibar2}), we only need to consider the coefficients of $q^{5n+2}$ in $(q^4;q^4)_\infty(q;q)_\infty^3$ and this was done in the latter part of the proof of Theorem \ref{congmod5} above.

\end{document}